
\def\LaTeX{%
  \let\Begin\begin
  \let\End\end
  \let\salta\relax
  \let\finqui\relax
  \let\futuro\relax}

\def\UK{\def\our{our}\let\sz s}
\def\USA{\def\our{or}\let\sz z}



\LaTeX

\USA


\salta


\documentclass[a4paper,12pt]{article}

\setlength{\textheight}{24cm}
\setlength{\textwidth}{16cm}
\setlength{\oddsidemargin}{2mm}
\setlength{\evensidemargin}{2mm}
\setlength{\topmargin}{-15mm}
\parskip2mm


\usepackage{amsmath}
\usepackage{amsthm}
\usepackage{amssymb}
\usepackage[mathcal]{euscript}
%
%
\usepackage{cite}
\usepackage[usenames,dvipsnames]{color}
%
%
%
\def\gianni{\color{black}}            
\def\betti{\color{black}} 
\def\juerg{\color{black}}              
\def\pier{\color{black}}
\def\gil{\color{black}}
\def\bet{\color{black}}
%
%
\let\gianni\relax
\let\betti\relax
\let\juerg\relax


\bibliographystyle{plain}


%

\finqui

\def\Beq{\Begin{equation}}
\def\Eeq{\End{equation}}
\def\Bsist{\Begin{eqnarray}}
\def\Esist{\End{eqnarray}}

\def\Bthm{\Begin{theorem}}
\def\Ethm{\End{theorem}}

\def\Brem{\Begin{remark}\rm}
\def\Erem{\End{remark}}

\def\Bcenter{\Begin{center}}
\def\Ecenter{\End{center}}
\let\non\nonumber




\def\step #1 \par{\medskip\noindent{\bf #1.}\quad}

\def\Step #1 \par{\bigskip\leftline{\bf #1}\nobreak\medskip}


\def\Lip{Lip\-schitz}
\def\holder{H\"older}
\def\aand{\quad\hbox{and}\quad}

\def\lhs{left-hand side}
\def\rhs{right-hand side}

\def\wk{well-known}


\def\generaliz{generali\sz}

\def\organiz{organi\sz}

\def\bhv{behavi\our}


\def\multibold #1{\def\arg{#1}%
  \ifx\arg\pto \let\next\relax
  \else
  \def\next{\expandafter
    \def\csname #1#1#1\endcsname{{\bf #1}}%
    \multibold}%
  \fi \next}

\def\pto{.}

\def\multical #1{\def\arg{#1}%
  \ifx\arg\pto \let\next\relax
  \else
  \def\next{\expandafter
    \def\csname cal#1\endcsname{{\cal #1}}%
    \multical}%
  \fi \next}


\def\multimathop #1 {\def\arg{#1}%
  \ifx\arg\pto \let\next\relax
  \else
  \def\next{\expandafter
    \def\csname #1\endcsname{\mathop{\rm #1}\nolimits}%
    \multimathop}%
  \fi \next}

\multibold
qwertyuiopasdfghjklzxcvbnmQWERTYUIOPASDFGHJKLZXCVBNM.

\multical
QWERTYUIOPASDFGHJKLZXCVBNM.

\multimathop
dist div dom meas sign supp .


\def\accorpa #1#2{\eqref{#1}--\eqref{#2}}
\def\Accorpa #1#2 #3 {\gdef #1{\eqref{#2}--\eqref{#3}}%
  \wlog{}\wlog{\string #1 -> #2 - #3}\wlog{}}


\def\separa{\noalign{\allowbreak}}

\def\todx{\mathrel{\scriptstyle\searrow}}

\def\QED{\hfill $\square$}

\def\graffe #1{\mathopen\{#1\mathclose\}}

\def\<#1>{\mathopen\langle #1\mathclose\rangle}
\def\norma #1{\mathopen \| #1\mathclose \|}
\def\normaX #1{{\mathopen \| #1\mathclose \|}_X}

\def\normaV #1{\norma{#1}_V}
\def\normaH #1{\norma{#1}_H}
\def\normaW #1{\norma{#1}_W}
\def\normaVp #1{\norma{#1}_*}

\def\iot {\int_0^t}
\def\ioT {\int_0^T}
\def\iO{\int_\Omega}
\def\intQt{\int_{Q_t}}

\def\dt{\partial_t}
\def\dn{\partial_\nu}

\def\cpto{\,\cdot\,}

\def\checkmmode #1{\relax\ifmmode\hbox{#1}\else{#1}\fi}

\def\aeQ{\checkmmode{a.e.\ in~$Q$}}
\def\aet{\checkmmode{a.e.\ in~$(0,T)$}}

\def\aat{\checkmmode{for a.a.~$t\in(0,T)$}}


\def\erre{{\mathbb{R}}}




\def\genspazio #1#2#3#4#5{#1^{#2}(#5,#4;#3)}
\def\spazio #1#2#3{\genspazio {#1}{#2}{#3}T0}

\def\L {\spazio L}
\def\H {\spazio H}

\def\C #1#2{C^{#1}([0,T];#2)}

\def\Vp{V^*}


\def\Lx #1{L^{#1}(\Omega)}
\def\Hx #1{H^{#1}(\Omega)}

\def\Luno{\Lx 1}
\def\Ldue{\Lx 2}

\def\Huno{\Hx 1}
\def\Hdue{\Hx 2}


\def\LQ #1{L^{#1}(Q)}


\let\a\alpha
\let\b\beta
\let\eps\varepsilon
\let\phi\varphi

\let\TeXchi\chi                         
\newbox\chibox
\setbox0 \hbox{\mathsurround0pt $\TeXchi$}
\setbox\chibox \hbox{\raise\dp0 \box 0 }
\def\chi{\copy\chibox}

\let\A\calA
\def\Am{\A^{-1}}
\def\Pi{\widehat\pi}


\def\hatC{\widehat C}


\def\muz{\mu_0}
\def\phiz{\phi_0}
\def\sigmaz{\sigma_0}
\def\muab{\mu_{\a,\b}}
\def\phiab{\phi_{\a,\b}}
\def\sigmaab{\sigma_{\a,\b}}
\def\xiab{\xi_{\a,\b}}
\def\Ra{R_{\a,\b}}

\def\mul{\overline\mu}
\def\phil{\overline\phi}
\def\sigmal{\overline\sigma}

\def\Rl{\overline R}

\def\psia{\psi_{\a,\b}}
\def\psil{\overline\psi}

\def\ustar{u^*}
\def\vstar{v^*}

\def\Rp{{\bet S}}
\def\mueps{\mu_\eps}
\def\phieps{\phi_\eps}
\def\xieps{\xi_\eps}
\def\sigmaeps{\sigma_\eps}
\def\Reps{R_\eps}
\def\muzeps{\mu_{0,\eps}}
\def\sigmazeps{\sigma_{0,\eps}}



\def\Beta{\widehat{\vphantom t\smash B\mskip2mu}\mskip-1mu}
\def\betaz{B^\circ}


\Begin{document}


\title{{{\betti Vanishing viscosities} and error estimate\\
for a Cahn--Hilliard type phase field system\\
related to tumor growth}%
\footnote{{\bf Acknowledgment.}\quad\rm  {\betti The financial support of the FP7-IDEAS-ERC-StG \#256872
(EntroPhase) is gratefully acknowledged by the authors. The present paper 
also benefits from the support of the MIUR-PRIN Grant 2010A2TFX2 ``Calculus of Variations'' for PC and GG, and the GNAMPA (Gruppo Nazionale per l'Analisi Matematica, la Probabilit\`a e le loro Applicazioni) of INdAM (Istituto Nazionale di Alta Matematica) for PC, GG and ER.}}}

\author{}
\date{}
\maketitle
\Bcenter
\vskip-2cm
{\large\sc Pierluigi Colli$^{(1)}$}\\
{\normalsize e-mail: {\tt pierluigi.colli@unipv.it}}\\[.25cm]
{\large\sc Gianni Gilardi$^{(1)}$}\\
{\normalsize e-mail: {\tt gianni.gilardi@unipv.it}}\\[.25cm]
{\large\sc Elisabetta Rocca$^{(2),\,(3)}$}\\
{\normalsize e-mail: {\tt {\betti rocca@wias-berlin.de and elisabetta.rocca@unimi.it}}}\\[.25cm]
{\large\sc J\"urgen Sprekels$^{(2),\,(4)}$}\\
{\normalsize e-mail: {\tt sprekels@wias-berlin.de}}\\[.45cm]
$^{(1)}$
{\small Dipartimento di Matematica ``F. Casorati'', Universit\`a di Pavia}\\
{\small via Ferrata 1, 27100 Pavia, Italy}\\[.2cm]
$^{(2)}$
{\small Weierstrass Institute for Applied Analysis and Stochastics}\\
{\small Mohrenstrasse 39, 10117 Berlin, Germany}\\[2mm]
$^{(3)}$
{\small Dipartimento di Matematica ``F. Enriques'', Universit\`a di Milano}\\
{\small Via Saldini 50, 20133 Milano, Italy}\\[2mm]
$^{(4)}$
{\small Department of Mathematics}\\
{\small Humboldt-Universit\"at zu Berlin}\\
{\small {\betti Unter den Linden 6, 10099} Berlin, Germany}\\
[1cm]
\Ecenter

\Begin{abstract}
{\betti In this paper  we  perform an asymptotic analysis for two different vanishing viscosity coefficients} {\juerg occurring in a phase field system of Cahn--Hilliard type that was recently} {\betti introduced in order to approximate a tumor growth model. In particular, we extend some recent results obtained in \cite{CGH}, letting the two positive viscosity parameters tend to zero independently from each other and weakening the conditions on the initial data in such a way {\juerg as} to maintain the nonlinearities of the PDE system as general as possible. Finally, under proper growth conditions on the {\pier \emph{interaction potential}}, we prove an error estimate leading also to the uniqueness result for the limit system. }
\vskip3mm

\noindent {\bf Key words:} {\betti tumor growth, Cahn--Hilliard {\pier system}, reaction-diffusion equation, asymptotic analysis, error estimates.}
\vskip3mm
\noindent {\bf AMS (MOS) Subject Classification:} {\betti 35Q92, 92C17, 35K35, 35K57, 78M35, 35B20, 65N15, 35R35.} 
\End{abstract}

\salta

\pagestyle{myheadings}
\newcommand\testopari{\sc Colli \ --- \ Gilardi \ --- \ Rocca \ --- \ Sprekels}
\newcommand\testodispari{\sc Asymptotic analysis for a Cahn--Hilliard type phase field system}
\markboth{\testodispari}{\testopari}

\finqui


\section{Introduction}
\label{Intro}
\setcounter{equation}{0}

{\betti In this paper we study}
the {\juerg system} of partial differential equations
\Bsist
  && \a \dt\mu + \dt\phi - \Delta\mu
  = p(\phi) (\sigma - \gamma\mu)
  \label{Iprima}
  \\
  && \mu = \b\dt\phi - \Delta\phi + F'(\phi)
  \label{Iseconda}
  \\
  && \dt\sigma - \Delta\sigma
  = - p(\phi) (\sigma - \gamma\mu),
  \label{Iterza}
\Esist
together with the boundary and initial conditions
\Bsist
  && \dn\mu = \dn\phi = \dn\sigma = 0 
  \label{Ibc}
  \\
  && \mu(0) = \muz, \quad
  \phi(0) = \phiz
  \aand
  \sigma(0) = \sigmaz \,.
  \label{Icauchy}
\Esist
\Accorpa\Ipbl Iprima Icauchy
Each of the partial differential equations \accorpa{Iprima}{Iterza}
is meant to hold in a three-dimensional bounded domain $\Omega${\pier ,}
endowed with a smooth boundary~$\Gamma${\pier ,} and for every positive time,
and $\dn$ in \eqref{Ibc} stands for the {\juerg outward} normal derivative on~$\Gamma$.
Moreover, $\a$~and $\b$ are nonnegative parameters, strictly positive in principle,
while $\gamma$ is a strictly positive constant.
Furthermore, $p$~is a nonnegative function, and {\betti $F$ is a nonnegative potential}.
Finally, $\muz$, ${\pier \phiz}$ and $\sigmaz$ are given initial data defined in~$\Omega$.

{\betti The physical context of this paper is that of {\juerg tumor growth} dynamics. This topic {\juerg has in recent years become of big interest} in applied mathematics, especially  {\juerg after continuum models were} developed (cf., e.g., \cite{CL,Letal}). 
{\juerg The fact that multiple constituents interact with each other made it necessary} to consider diffuse interface models based on continuum mixture theory (cf., e.g., \cite{CBCB,CLLW,Fetal,HPZO,OHP,WLFC}). These models consist of a Cahn--Hilliard type equation (in general with transport) containing reaction terms {\juerg that depend on the nutrient concentration (e.g.~oxygen) and in turn obey} an advection-reaction-diffusion equation. 
Even {\juerg though} numerical simulations of these models have already been carried out in several papers (cf., e.g., \cite[Chapter~8]{CL} and references therein), the rigorous mathematical analysis of the resulting PDEs systems is still very poor. To our knowledge, the first results are related to the so-called Cahn--Hilliard--Hele--Shaw system (cf., e.g., \cite{LTZ,WZ}) {\juerg in which} the nutrient is neglected, while two very recent contributions \cite{CGH} and \cite{FGR}  (cf. also  \cite{HKNV}, where formal studies on the corresponding sharp interface limits are performed) deal with a model recently proposed in \cite{HZO} (or approximations {\juerg thereof}, see also \cite{WZZ}), where the velocities are set to zero and the state variables are the tumor fraction $\varphi$ and the nutrient-rich concentration $\sigma$.  We can set $\phi\simeq 1$ in the tumorous phase and $\phi\simeq -1$ in the healthy cell phase, while $\sigma$ {\juerg typically satisfies} $\sigma\simeq 1$ in a nutrient-rich extracellular water phase 
and $\sigma\simeq 0$ in a nutrient-poor extracellular water phase. Moreover, 
the third unknown $\mu$ is the related chemical potential, specified by \eqref{Iseconda} 
as in the case of the viscous Cahn--Hilliard or Cahn--Hilliard equation, depending on whether 
$\b>0 $ or $\b=0$ (see \cite{CH, EllSt, EllZh}). {\juerg In addition, in \cite{CGH} the PDE system \Ipbl\ was studied for the very particular case that $\a=\b$, and 
the asymptotic analysis as the coefficient $\a=\b$ tends to zero was performed, yielding 
the convergence of subsequences to weak solutions of the limit problem; moreover, in 
\cite{FGR} the existence of weak solutions, as well as uniqueness and existence of attractors, was proved directly} for the limit system where $\a=\b=0$ (cf.~also the following comments in this Introduction).}

In the case $\a=0$, 
the sub-system \accorpa{Iprima}{Iseconda} becomes of viscous or pure Cahn--Hilliard type,
depending on whether $\b>0$ or $\b=0$. 
On the other hand, {\juerg in the case} $\a>0$ the presence of the term $\alpha\dt\mu$ in \eqref{Iprima} 
gives a parabolic structure to equation~\eqref{Iprima} with respect to $\mu$. 

We remark that the original model deals with functions
$F$~and $p$ that are precisely related to each other. 
Namely, we~have
\Beq
  p(u) = 2p_0 \sqrt{F(u)}
  \quad \hbox{if $|u|\leq 1$}
  \aand
  p(u) = 0
  \quad \hbox{otherwise,}
  \label{relazWp}
\Eeq
where $p_0$ is a positive constant and $F(u)$ is the classical
Cahn--Hilliard double-well free energy density. 
However, this relation is useless in many aspects of the mathematical study. 
Moreover, one can allow $F$ to be even a singular potential.

{\juerg As mentioned above, \cite{CGH}~just deals with the case $\a=\b$ for the mathematical study, although the constants $\a$ and $\b$ have a different meaning.
In that paper, the existence of a unique solution to the system~\Ipbl\ was proved 
under very general conditions on $p$ and~$F$, and, in the same framework, the long-time \bhv\ of the solution was discussed. In addition,
in a more restricted setting for the double-well potential~$F$, 
\cite{CGH}~investigated the asymptotic \bhv\ of the problem as the coefficient $\a=\b$ tends to zero, finding} the convergence of subsequences to weak solutions of the limit problem.
Moreover, under a smoothness condition on the initial values{\pier ,}
uniqueness for the limit problem {\juerg was} proved and, consequently, also the convergence of the entire family.
It must be pointed out that a uniqueness result {\juerg was} proved in \cite{FGR} under weaker assumptions.

In the present paper, we first extend some of the results of \cite{CGH}.
Namely, we let the positive parameters $\a$ and $\b$ be independent from each other,
and we weaken the assumptions on the initial data {\juerg while}
keeping the potential as general as possible.
At the same time, we establish a general a priori estimate
that is uniform with respect to the parameters $\a$ and~$\b$.
This is the starting point of possible asympotic analyses
with respect to {\juerg these} parameters.
Then, we confine ourselves to a class of regular potentials.
In this framework, we state a convergence result as both $\a$ and $\b$ tend to zero independently,
and we prove an error estimate in terms of $\a$ and $\b$
for the difference of the solution to \Ipbl\ and the one of the limit problem.
The case of just one of the parameters tending to zero is the subject of a work in progress.

{\pier Let us express our belief {\bet that the results of} the present paper are general and interesting enough so that methods and estimates could be extended to other situations. In particular,
in case of the trivial choice $p\equiv 0$ (admitted by our assumption \eqref{hpconst})
our system \eqref{Iprima}--\eqref{Iterza} decouples and \eqref{Iprima}--\eqref{Iseconda}
reduces to a well-known phase field system of Caginalp tye which can be seen as a 
(doubly) viscous approximation of the Cahn--Hilliard system 
obtained at the limit as $\alpha $ and $\beta$ go to zero. 
To this concern, let us quote the papers \cite{DHK, DHL, Ros1, Ros2},  where 
different investigations on this kind of viscous approximations of Cahn--Hilliard system 
are performed, and point out that the results contained in \cite{Ros1} are here generalized 
and somehow improved.} 

Our paper is \organiz ed as follows. 
In the next section, we will state the assumptions and our results on the mathematical problem. 
In Section~\ref{EXTENSIONS}, we will prove the extensions mentioned above.
The last section is devoted to the asymptotic analysis and the error estimate.
In the {\juerg remainder} of the paper, we take $\gamma=1$, without loss of generality.


\section{Statement of the problem and results}
\label{STATEMENT}
\setcounter{equation}{0}

In this section, we make precise assumptions and state our results.
As in the Introduction, $\Omega\subset\erre^3$
{\juerg denotes} the domain where the evolution takes place and $\Gamma$ {\pier is}
{\juerg its} boundary.
We assume $\Omega$ to be open, bounded, and connected, and $\Gamma$ to be smooth.
Moreover, the symbol $\dn$ denotes the {\juerg outward} normal derivative on~$\Gamma$.
Given a final time~$T$, we~set
\Beq
  Q := \Omega\times(0,T)
  \aand
  \Sigma := \Gamma \times (0,T) .
  \label{defQ}
\Eeq
Moreover, we set for brevity
\Beq
  V := \Huno,
  \quad H := \Ldue ,
  \aand
  W := \graffe{v\in\Hdue:\ \dn v = 0 \ \hbox{on $\Gamma$}},
  \label{defspazi}
\Eeq
and endow {\pier these spaces with their} standard norms.
For the norm in a generic Banach space~$X$ (or~a power of~it), we use the symbol~$\norma\cpto_X$
with the following exceptions: we simply write $\norma\cpto_p$ and $\normaVp\cpto$
if $X=\Lx p$ or $X=\LQ p$ for $p\in[1,+\infty]$ and $X=\Vp$, the dual space of~$V$, respectively.
Finally, it is understood that {\juerg $H$ is embedded in $\Vp$ in the usual way, i.e., such}
 that $\<u,v>=\iO u\,v\,$ for every $u\in H$ and $v\in V$,
where $\<\cpto,\cpto>$ stands for the duality pairing between $\Vp$ and~$V$.

\bigskip

As far as the structure of the system is concerned, we are given
two constants $\a$ and $\b$ and three functions $p$, $\Beta$ and $\Pi$
satisfying the conditions listed below
\Bsist
  && \a ,\, \b \in (0,1) 
  \label{hpconst}
  \\
  && \hbox{$p:\erre\to{\betti [0,+\infty)}$ is bounded and \Lip\ continuous}
  \label{hpp}
  \\
  && \hbox{$\Beta:\erre\to[0,+\infty]$ is convex, proper, lower
      semicontinuous}
    \qquad
  \label{hpBeta}
  \\
  && \hbox{$\Pi\in C^1(\erre)$ is nonnegative, and
      $\pi:=\Pi\,'$ is \Lip\ continuous}.
  \label{hpPi}
\Esist
We also define the potential $F:\erre\to[0,+\infty]$
and the graph $B$ in $\erre\times\erre$ by
\Beq
  F := \Beta + \Pi 
  \aand
  B := \partial\Beta .
  \label{defWbeta}
\Eeq
\Accorpa\HPstruttura hpconst defWbeta
We notice that {\juerg  if $F$ is a $C^2$ function then}
our assumptions imply that $F''$ is bounded from below.
We also remark that $B$ is maximal monotone.
In the {\juerg following}, we write $D(\Beta)$ and $D(B)$ for the effective domains
of~$\Beta$ and~$B$, respectively, and we use the same symbol $B$
for the maximal monotone operators induced on $L^2$ spaces.

\Brem
\label{CompatWp}
We notice that, among many others,
the most important and typical examples of potentials fit our assumptions.
Namely, we can take as $F$ the classical double-well potential
and as the logarithmic {\juerg potential, which are} defined~by
\Bsist
  \hskip-.8cm && F_{cl}(r) := {\textstyle \frac 14} (r^2-1)^2
  = {\textstyle \frac 14} ((r^2-1)^+)^2 + {\textstyle \frac 14} ((1-r^2)^+)^2
  \quad \hbox{for $r\in\erre$}
  \label{clW}
  \\
  \hskip-.8cm && F_{log}(r) := (1-r)\ln(1-r) + (1+r)\ln(1+r) + \kappa (1 - r^2 )^+
  \quad \hbox{for $|r|<1$},
  \qquad
  \label{logW}
\Esist
where the decomposition $F=\Beta+\Pi$ as in \eqref{defWbeta} is written {\juerg explicitly}.
In \eqref{logW}, $\kappa$~is a positive constant
which{\juerg, depending on its value, does or does not provide a double well},
and the definition of the logarithmic part of $F_{log}$
is extended by continuity {\juerg to} $\pm1$ and by $+\infty$ outside~$[-1,1]$.
Moreover, another possible choice is 
\Beq
  F(r) := I(r) + ((1 - r^2)^+)^2
  \quad \hbox{for $r\in\erre$},
  \label{irrW}
\Eeq
where $I$ is the indicator function of $[-1,1]$,
{\juerg which takes} the value $0$ in $[-1,1] $ and $+\infty$ elsewhere.
 {\juerg For such an irregular potential, the associated}
subdifferential is multi-valued, and the precise statement of problem \Ipbl\
has to introduce a selection $\xi$ of~$B(u)$.
\Erem

As far as the initial data of our problem are concerned,
we assume that
\Beq
  {\betti \sqrt{\alpha}} \, \muz \, , \, \sigmaz \in H , \quad
  \phiz \in V,
  \aand
  F(\phiz) \in \Luno ,
  \label{hpdati}
\Eeq
while the regularity properties {\juerg postulated} for the solution
are the following:
\Bsist
  && \mu , \sigma \in \H1\Vp \cap \L2V 
  \label{regmusigma}
  \\
  && \phi \in \H1H \cap \L2W 
  \label{regphi}
  \\
  && \xi \in \L2H ,
  \aand
  \xi \in B(u)
  \quad \aeQ .
  \label{regxi}
\Esist
\Accorpa\Regsoluz regmusigma regxi
We notice that \accorpa{regmusigma}{regphi} imply {\juerg that}
$\mu,\sigma\in\C0H$ and $\phi\in\C0V$.
At this point, we consider the problem 
of finding a quadruplet $(\mu,\phi,\sigma,\xi)$
with the above regularity 
in order that $(\mu,\phi,\sigma,\xi)$ and the related function
\Beq
  R = p(\phi) (\sigma - \mu)
  \label{defR}
\Eeq
satisfy the system
\Bsist
  \hskip-1cm
  && \a \< \dt\mu , v > + \iO \dt\phi \, v + \iO \nabla\mu \cdot \nabla v
  = \iO R v 
  \non
  \\
  && \quad \hbox{for every $v\in V$, \ \aet} 
  \label{prima}
  \\ [3pt]
  && \mu = \b\dt\phi - \Delta\phi + \xi + \pi(\phi)
  \aand
  \xi \in B(\phi)
  \quad \aeQ
  \label{seconda}
  \\
  && \< \dt\sigma , v > + \iO \nabla\sigma \cdot \nabla v
  = - \iO R v 
  \non
  \\
  && \quad \hbox{for every $v\in V$, \ \aet} 
  \label{terza}
  \\ [3pt]
  && \mu(0) = \muz, \quad
  \phi(0) = \phiz
  \aand
  \sigma(0) = \sigmaz \,.
  \label{cauchy}
\Esist
\Accorpa\Pbl defR cauchy
\Accorpa\Tuttopbl regmusigma cauchy
This is a weak formulation of the boundary value problem \Ipbl\ described in the Introduction.
The homogeneous Neumann boundary condition for $\phi$ 
is contained in~\eqref{regphi} (see \eqref{defspazi} for the definition of~$W$),
while the analogous ones for $\mu$ and $\sigma$
are meant in a \generaliz ed sense through the variational equations \eqref{prima} and~\eqref{terza}.
We notice once and for all that {\pier the addition} of \eqref{prima} and~\eqref{terza} yields
\Beq
  \< \dt \bigl( \a \mu + \phi + \sigma \bigr) , v >
  + \iO \nabla (\mu + \sigma) \cdot \nabla v
  = 0
  \label{unomenotre}
\Eeq
for every $v\in V$, \aet.
We also set for convenience
\Beq
  \Rp = \sqrt{p(\phi)} \, (\sigma - \mu) \,.
  \label{defRp}
\Eeq

Our first results deal with {\juerg the well-posedness of} the above problem
and general a priori estimates.
Namely, we have:

\Bthm
\label{Wellposedness}
Assume \HPstruttura\ and \eqref{hpdati}.
Then, for every $\a,\b\in(0,1)$, there exists a unique quadruplet $(\mu,\phi,\sigma,\xi)$
satisfying \Regsoluz\ and solving problem \Pbl.
\Ethm

\Bthm
\label{GenEst}
Assume \HPstruttura\ and \eqref{hpdati}.
Then, for some constant $\hatC$ that depends only on 
$\Omega$, $T$ and the {\betti shapes of~$\pi$ and $p$},
the following is {\juerg true: for} every $\a,\b\in(0,1)$,
the solution $(\mu,\phi,\sigma,\xi)$
to problem \Pbl\ with the regularity specified by~\Regsoluz\ satisfies
\Bsist
  && \a^{1/2} \norma\mu_{\L\infty H}
  + \norma{\nabla\mu}_{\L2H}
  \non
  \\
  && \quad {}
  + \b^{1/2} \norma{\dt\phi}_{\L2H}
  + \norma\phi_{\L\infty V}
  + \norma{F(\phi)}_{\L\infty\Luno}^{1/2}
  \non
  \\
  && \quad {}
  + \norma\sigma_{\H1\Vp\cap\L\infty H\cap\L2V}
  + \norma\Rp_{\L2H}
  + \norma R_{\L2H}
  \non
  \\
  && \quad {}
  + \norma{\dt(\a\mu+\phi)}_{\L2\Vp}
  \non
  \\
  && \leq \hatC \, \bigl(
    \a^{1/2} \normaH\muz + \normaV\phiz + \norma{F(\phiz)}_{\Luno}^{1/2} + \normaH\sigmaz 
  \bigr)
  \label{genest}
\Esist
as well~as
\Bsist
  && \norma\mu_{\L2V}
  + \norma\phi_{\L2W}
  + \norma\xi_{\L2H}
  \non
  \\
  && \leq \hatC \, \bigl(
    \a^{1/2} \normaH\muz + \normaV\phiz + \norma{F(\phiz)}_{\Luno}^{1/2} + \normaH\sigmaz + \norma\mu_{\L2H} + 1
  \bigr) .
  \qquad
  \label{genestbis}
\Esist
\Ethm

Thus, a uniform estimate for the \lhs\ of \eqref{genest}
holds in terms of the norms of the initial data related to~\eqref{hpdati},
while an estimate for the \lhs\ of \eqref{genestbis}
follows whenever a bound for $\norma\mu_{\L2H}$ {\juerg has been} proved.

\Brem
\label{WellposCGH}
We note that Theorems~\ref{Wellposedness} and~\ref{GenEst}
improve the results of~\cite{CGH},
since the stronger assumption made there,
\Beq
  \muz, \phiz, \sigmaz \in V
  \aand
  F(\phiz) \in \Luno,
  \label{hpdatiCGH}
\Eeq
is now replaced by~\eqref{hpdati}
(and also since just the case $\a=\b$ is dealt with in~\cite{CGH}).
\Erem

Our next results regard the asymptotic analysis
as the coefficients $\a$ and $\b$ tend to zero, independently.
To this end, we restrict ourselves
to a particular class of potentials.
Namely, we also assume~that
{\betti
\Beq
  \gianni{%
  D(\Beta) = \erre
  \aand
  |\betaz(r)|\leq C\left(\widehat B(r)+1\right)
  \quad \hbox{for every $r\in\erre$,}
  }
  \label{hpBbase}
\Eeq
where $\betaz$ is the element of $B$ with minimal norm and $C$ is a given positive constant. Let us note {\juerg that,
 for example, all polynomially growing potentials, as well as exponential functions,} comply with our assumption \eqref{hpBbase}}.  {\pier Let us point out that \eqref{hpBbase} implies (actually, it is equivalent to)  the condition
\Beq
  D(\Beta) = \erre, 
  \quad 
  |s|\leq C\left(\widehat B(r)+1\right)
  \quad \hbox{for all $r\in\erre$, $s\in B (r)$}
  \label{hp-pier}
\Eeq
for the same constant $C$, as checked precisely in the next remark.} 
 

{\gil
\Brem
\label{GenCondCrescita}
In fact, a~similar equivalence holds 
for a more general growth condition
and in the general setting of Hilbert spaces,
as we show at once.
If ${\pier X}$~is a Hilbert space, $\Beta : {\pier X}\to[0,+\infty)$ is convex and~l.s.c.\
(thus continuous since it is everywhere defined),
$B :=\partial\Beta$ and, for every $u\in {\pier X}$, 
$\betaz(u)$~is the element of $B (u)$ having {\bet minimal} norm,
the assumption
\Beq
  \normaX{\betaz(u)} \leq {\pier \Psi}  (\Beta(u))
  \quad \hbox{for every $u\in {\pier X}$,}
  \non
\Eeq
where ${\pier \Psi} :[0,+\infty)\to[0,+\infty)$ is continuous,
implies
\Beq
  \normaX{{\pier\zeta}} \leq {\pier \Psi} (\Beta(u))
  \quad \hbox{for every $u\in {\pier X}$ and every ${\pier\zeta}\in B(u) $}.
  \non
\Eeq
Indeed, for arbitrary $u\in {\pier X}$, ${\pier\zeta}\in B(u)$ and $\eps>0$, we have
\Beq
  \bigl( \betaz(u+\eps{\pier\zeta}) - {\pier\zeta} , (u+\eps{\pier\zeta}) - u \bigr) \geq 0 ,
  \quad \hbox{whence} \quad
  \normaX{{\pier\zeta}} \leq \normaX{\betaz(u+\eps{\pier\zeta})}.
\Eeq
By applying our assumption to $u+\eps{\pier\zeta}$, we deduce that
\Beq
  \normaX{{\pier\zeta}} \leq {\pier \Psi} (\Beta(u+\eps{\pier\zeta})).
  \non
\Eeq
By taking $\eps\to0$ and owing to the continuity of ${\pier \Psi} \circ\Beta$, we conclude.
\Erem
}

Now we are ready to state  our result on asymptotics.
\Bthm
\label{Asymptotics}
Assume \HPstruttura\ and {\gianni \eqref{hpBbase} on the structure}
and \eqref{hpdati} on the initial data.
Moreover, let {\betti $(\muab,\phiab,\sigmaab,\xi_{\alpha,\beta})$} be the unique solution to problem \Pbl\
given by Theorem~\ref{Wellposedness}.
Then, we have {\betti that there exists a quadruplet $(\mu,\phi,\sigma,\xi)$ such that}
\Bsist
  & \muab \to \mu
  & \hbox{weakly in $\L2V$}
  \label{convmu}
  \\
  & \phiab \to \phi
  & \hbox{weakly star in $\L\infty V\cap\L2W$} \qquad
  \label{convphi}
  \\
  & \sigmaab \to \sigma
  & \hbox{weakly in $\H1\Vp\cap\L2V$}
  \label{convsigma}
  \\
  & \dt (\a\muab+\phiab) \to \dt\phi
  & \hbox{weakly in $\L2\Vp$}
  \label{convdt}
   \\
  & {\pier \xiab \to \xi}
  &{\pier \hbox{weakly in $\L2H$}}
  \label{pier1}
\Esist
at least for a subsequence.
Moreover, every limiting {\betti quadruplet $(\mu,\phi,\sigma,\xi)$}
satisfies
\Bsist
  & \< \dt\phi , v > + \iO \nabla\mu \cdot \nabla v
  = \iO R \, v
  & \quad \hbox{$\forall\,v\in V$, \ \aet}
  \label{primalim}
  \\
  & \mu = - \Delta\phi {\betti {} + \xi+{\pier \pi } (\phi), \quad \xi\in B(\varphi)}
  & \quad \aeQ
  \label{secondalim}
  \\
  & \< \dt\sigma , v > + \iO \nabla\sigma \cdot \nabla v
  = - \iO R \, v
  & \quad \hbox{$\forall\,v\in V$, \ \aet}
  \label{terzalim}
  \\
  & \phi(0) = \phiz
  \aand
  \sigma(0) = \sigmaz
  & \quad \hbox{in $\Omega$} 
  \label{cauchylim}
\Esist
where $R$ is defined by~\eqref{defR}, accordingly.
\Ethm

The above result \generaliz es the analogous \cite[Thm.~2.6]{CGH}
as far as the assumptions on the initial data are concerned
(and also since just the case $\a=\b$ was considered there).
Moreover, {\betti in \cite[Thm.~2.6]{CGH}},
even uniqueness for {\juerg the solution to the limit problem was} proved.
However, also for this point, 
stronger conditions on the initial data are assumed
in order that the solution to the limit problem is rather smooth.
Here, we can consider
the natural regularity requirements,~i.e.,
\Bsist
  && \mu \in \L2V
  \label{regmulim}
  \\
  && \phi \in \H1\Vp {\pier {} \cap \L\infty V } \cap \L2W 
  \label{regphilim}
  \\
  && \sigma \in \H1\Vp \cap \L2V \subset \C0H \,.
  \label{regsigmalim}
\Esist
\Accorpa\Regsoluzlim regmulim regsigmalim
For uniqueness in this framework, we can {\pier quote} the even more general result~\cite[Thm.~2]{FGR}.
However, uniqueness also follows from the error estimate we present at once
(see the forthcoming Remark~\ref{Uniqueness} for details). {\betti In order to state our last result we need to {\pier reinforce} the 
assumptions we made on the potential $F${\pier ;} namely, we assume
\Bsist
   D(\Beta) = \erre
  \aand
  \hbox{$F=\Beta+\Pi$ \enskip is a $C^2$ function on $\erre$}
  \qquad\quad
  \label{ovunque}
  \\
  |F(r)| \leq C_0 (|r|^6 + 1), \quad
  |F'(r)| \leq C_1 (|r|^5 + 1), 
  \aand
  |F''(r)| \leq C_2 (|r|^4 + 1).
  \quad \hbox{}
  \label{classico}
\Esist
}%
{\juerg Although the third condition in \eqref{classico} implies the other two},
we have written all of them for convenience.
We also remark that the classical potential~\eqref{clW} fulfils such assumptions.
Furthermore, we notice that \eqref{classico} is slightly more general than
the analogous assumption made in \cite[Thm.~2.6]{CGH}. {\betti Finally, we can observe that the exponents 
in \eqref{classico} are related to the dimension of $\Omega$ and the {\juerg related} Sobolev embeddings}.
Here is our last result.

\Bthm
\label{Error}
Assume \HPstruttura\ and {\betti \accorpa{ovunque}{classico} {\gianni on the structure}
and \eqref{hpdati}} {\gianni on the initial data}.
Then, with the notation of Theorem~\ref{Asymptotics}, 
the estimate
\begin{align}
  \norma{\phiab-\phi}_{\L\infty\Vp\cap\L2V}
  {\pier{} + \norma{\muab-\mu}_{\L2\Vp}}
  \non
  \\
  + \norma{\sigmaab-\sigma}_{\L\infty\Vp\cap\L2H}
  \leq C \, \bigl( \a^{1/2} + \b^{1/2} \bigr)
  \label{error}
\end{align}
holds true with a constant $C$ that depends only on~$\Omega$, $T$, 
the structure of the system,
and the norms of the initial data related to assumptions~\eqref{hpdati}, {\betti but not on $\alpha$ nor on $\beta$}.
\Ethm

The rest of the section is devoted to list some facts.
We {\juerg make repeated} use of the notation
\Beq
  Q_t := \Omega \times (0,t)
  \quad \hbox{for $t\in[0,T]$}
  \label{defQt}
\Eeq
and of \wk\ inequalities, namely, of 
the elementary Young inequality
\Beq
  ab \leq \delta a^2 + \frac 1{4\delta} \, b^2
  \quad \hbox{for every $a,b\geq 0$ and $\delta>0$}
  \label{young}
\Eeq
{\juerg as well as of H\"older's} inequality and its consequences.
Moreover, as $\,\Omega\,$ is bounded and smooth,
we can owe to the Poincar\'e and Sobolev type inequalities,
namely, 
\Bsist
  && \normaV v \leq C \Bigl( \normaH{\nabla v} + \bigl| \textstyle\iO v \bigr| \Bigr)
  \quad \hbox{for every $v\in V$} 
  \label{poincare}
  \\
  && V \subset \Lx q
  \aand
  \norma v_q \leq C \normaV v
  \quad \hbox{for every $v\in V$ and $1\leq q\leq 6$}
  \label{sobolev}
  \\
  && \Lx q \subset \Vp
  \aand
  \normaVp v \leq C \norma v_q
  \quad \hbox{for every $v\in\Lx q$ and $q\geq 6/5$} \,.
  \qquad
  \label{dualsobolev}
\Esist
In \accorpa{poincare}{dualsobolev}, $C$ only depends on~$\Omega$.
Finally, we recall the interpolation inequality
\Beq
  \normaH v^2 \leq \normaV v \, \normaVp v
  \quad \hbox{for every $v\in V$,}
  \label{interpol}
\Eeq
which trivially follows from the identity 
$\normaH v^2=\<v,v>$ for every $v\in V$.


\section{Proofs of Theorems \ref{Wellposedness} and \ref{GenEst}}
\label{EXTENSIONS}
\setcounter{equation}{0}

We start proving Theorem~\ref{GenEst} in the following form:
\accorpa{genest}{genestbis}
hold for every $\a$ and $\b$ and every solution to problem \Pbl\
satisfying the regularity specified by~\Regsoluz.
We do not know anything about well-posedness yet, indeed.

\step
First a priori estimate

We test \eqref{prima} and \eqref{terza} by $\mu$ and~$\sigma$, respectively,
and integrate over~$(0,t)$, where $t\in(0,T)$ is arbitrary.
At the same time, we multiply \eqref{seconda} by $-\dt\phi$ and integrate over~$Q_t$.
Then, we add the resulting equalities to each other, {\betti obtaining}
\Bsist
  && \frac \a 2 \iO |\mu(t)|^2
  + \intQt \dt\phi \, \mu
  + \intQt |\nabla\mu|^2
  \non
  \\
  && \quad {}
  - \intQt \mu \, \dt\phi
  + \b \intQt |\dt\phi|^2
  + \frac 12 \iO |\nabla\phi(t)|^2
  + \iO F(\phi(t))
  \non
  \\
  && \quad {}
  + \frac 12 \iO |\sigma(t)|^2
  + \intQt |\nabla\sigma|^2 
  + \intQt R(\sigma-\mu)
  \non
  \\
  && = \frac \a 2 \iO |\muz|^2
  + \frac 12 \iO |\nabla\phiz|^2
  + \iO F(\phiz)
  + \frac 12 \iO |\sigmaz|^2 \,.
  \non
\Esist
Clearly, two terms cancel out.
Moreover, $F$~is nonnegative by assumptions \accorpa{hpBeta}{hpPi}.
Finally, we have $R(\sigma-\mu)=|\Rp|^2$
and $|R|\leq|\Rp|\sup\sqrt p$ \ \aeQ\ with the notation~\eqref{defRp}.
Therefore, {\pier with the help of \eqref{hpp}} we immediately deduce 
\Bsist
  && \a^{1/2} \norma\mu_{\L\infty H}
  + \norma{\nabla\mu}_{\L2H}
  \non
  \\
  && \quad {}
  + \b^{1/2} \norma{\dt\phi}_{\L2H}
  + \norma{\nabla\phi}_{\L\infty H}
  + \norma{F(\phi)}_{\L\infty\Luno}^{1/2}
  \non
  \\
  && \quad {}
  + \norma\sigma_{\L\infty H\cap\L2V}
  + \norma\Rp_{\L2H}
  + \norma R_{\L2H}
  \non
  \\
  && \leq C \, \bigl(
    \a^{1/2} \normaH\muz + \normaH{\nabla\phiz} + \norma{F(\phiz)}_{\Luno}^{1/2} + \normaH\sigmaz \bigr)
  \label{quasigenest}
\Esist
{\pier for some} {\gianni constant~$C$ that depends only on $p$}.
Thus, in order to prove~\eqref{genest},
we have to complete the full norm of $\phi$ 
and estimate the {\juerg terms} that are missing in~\eqref{quasigenest}.

\step
Second a priori estimate

We estimate the mean value of~$\phi$ by testing \eqref{unomenotre} by $v=1$.
We obtain, for every $t\in[0,T]$,
\Beq
  \iO \bigl( \a\mu(t) + \phi(t) + \sigma(t) \bigr)
  = \iO \bigl( \a\muz + \phiz + \sigmaz \bigr)
  \leq |\Omega|^{1/2} \normaH{\a\muz + \phiz + \sigmaz}
  \non
\Eeq
and deduce that (since $\alpha<1$)
\Beq
  \Bigl| \iO \phi(t) \Bigr|
  \leq  C_\Omega \bigl(
    \a^{1/2} \normaH\muz + \normaH\phiz + \normaH\sigmaz
    + \a^{1/2} \normaH{\mu(t)}
    + \normaH{\sigma(t)}
  \bigr), 
  \label{estmeanphi}
\Eeq
where $C_\Omega$ depends only on~$\Omega$.

\step
Third a priori estimate

We test \eqref{terza}, written at the time $t$, with $v(t)$,
where $v$ is arbitrary in $\L2V$.
Then we integrate over~$(0,T)$ with respect to~$t$ and obtain
\Beq
  \Bigl| \ioT \<\dt\sigma(t) , v(t) > \, dt \Bigr|
  \leq \bigl( \norma{\nabla\sigma}_{\L2H} + \norma R_{\L2H} \bigr) \norma v_{\L2V} \,.
  \non
\Eeq
This means that
\Beq
  \norma{\dt\sigma}_{\L2\Vp}
  \leq \norma{\nabla\sigma}_{\L2H} + \norma R_{\L2H} \,.
  \label{estdtsigma}
\Eeq

\step
Fourth a priori estimate

Similarly, we test \eqref{unomenotre}, written at the time $t$, by~$v(t)$,
where $v$ is arbitrary in $\L2V$.
We obtain
\Beq
  \Bigl| \ioT \<\dt(\a\mu+\phi)(t) , v(t) > \, dt \Bigr|
  \leq \Bigl| \ioT \<\dt\sigma(t) , v(t) > \, dt \Bigr|
  + \norma{\nabla(\mu+\sigma)}_{\L2H} \norma v_{\L2V}\,, 
  \non
\Eeq
whence immediately
\Beq
  \norma{\dt(\a\mu+\phi)}_{\L2\Vp}
  \leq \norma{\dt\sigma}_{\L2\Vp}
  + \norma{\nabla\mu}_{\L2H}
  + \norma{\nabla\sigma}_{\L2H} .
  \label{estdtaltro}
\Eeq

\step
First conclusion

We combine \accorpa{quasigenest}{estdtaltro}
with the Poincar\'e inequality~\eqref{poincare} applied to $\phi$
and immediately deduce~\eqref{genest}
with a constant $\hatC$ that depends only on {\betti $p$,} $\Omega$ and~$T$.

\step
Fifth a priori estimate and conclusion

By estimate~\eqref{genest} and the \Lip\ continuity of~$\pi$,
we deduce {\juerg that}
\Beq
  \norma{\pi(\phi)}_{\L2H} 
  \leq \hatC \, \bigl(
    \a^{1/2} \normaH\muz + \normaV\phiz + \norma{F(\phiz)}_{\Luno}^{1/2} + \normaH\sigmaz + 1
  \bigr),
  \label{estpi}
\Eeq
with the same $\hatC$, without loss of generality, provided that
we allow $\hatC$ to depend on~$\pi$ as well.
Now, we write \eqref{seconda} in the~form
\Beq
  - \Delta\phi + \xi
  = f := - \b\dt\phi - \pi(\phi) + \mu
  \non
\Eeq
and observe that \eqref{genest}{\pier ,  \eqref{estpi} and $\beta<1$} imply
\Beq
  \norma f_{\L2H}
  \leq \hatC \, \bigl(
    \a^{1/2} \normaH\muz + \normaV\phiz + \norma{F(\phiz)}_{\Luno}^{1/2} + \normaH\sigmaz + 1 + \norma\mu_{\L2H}
  \bigr),
  \non
\Eeq
with the same $\hatC$ once more, without loss of generality.
If $M$ denotes the \rhs\ of this inequality,
a~standard argument (formally multiply by $-\Delta\phi$)
shows that both $\Delta\phi$ and $\xi$ are bounded in $\LQ2$ by a multiple of~$M$.
Therefore, the same holds for $\norma\phi_{\L2W}$ by elliptic regularity.
Finally, the full norm $\norma\mu_{\L2V}$ is is equivalent
to the sum of $\norma{\nabla\mu}_{\L2H}$ and~$\norma\mu_{\L2H}$.
Thus, \eqref{genestbis} follows and the proof of Theorem~\ref{GenEst} {\betti is complete}.
{\juerg \hfill{\qed}}

\step
Proof of Theorem~\ref{Wellposedness}

As far as uniqueness is concerned, we can refer to the proof 
of the uniqueness part of~\cite[Thm.~2.2]{CGH}
since it holds under the present assumptions.
In order to prove the existence of a solution,
we approximate the data $\muz$ and $\sigmaz$
by functions $\muzeps$ and~$\sigmazeps$ satisfying
\Beq
  \muzeps , \sigmazeps \in V \quad \hbox{for $\eps>0$}, \qquad
  \muzeps \to \muz
  \aand
  \sigmazeps \to \sigmaz
  \quad \hbox{in $H$\quad as $\eps\todx0$}.
  \non
\Eeq
Then, for every $\eps>0$, the condition \eqref{hpdatiCGH} holds for the approximating data
so that the assumptions of \cite[Thm.~2.2]{CGH} are fulfilled.
Thus, {\pier the} problem \Pbl\ has a unique solution $(\mueps,\phieps,\sigmaeps,\xieps)$
with $\Reps$ defined by \eqref{defR} accordingly.
Moreover, such a solution must satisfy \accorpa{genest}{genestbis}
due to the above proof.
As $\a$ and $\b$ are fixed, such estimates
provide uniform boundedness with respect to $\eps$
even for $\mueps$ {\pier in $\L\infty H$ and $\dt\phieps$ in $\L2H$}.
Therefore, \eqref{genestbis} implies that $\mueps$, $\phieps$ and $\xieps$
are bounded in $\L2V$, $\H1H\cap\L2W$ and~$\L2H$, respectively.
Finally, the estimate for the time derivative of $\a\mueps+\phieps$
derived from \eqref{genest} and the estimate for $\dt\phieps$ mentioned before
imply that $\dt\mueps$ is bounded in~$\L2\Vp$.
Hence, we~have
\begin{align}
  & \mueps \to \mu
  \quad \hbox{weakly star in $\H1\Vp\cap\L2V$}
  \non
  \\
  & \phieps \to \phi 
  \quad \hbox{weakly in $\H1H\cap\L2W$}
  \non
  \\
  & \sigmaeps \to \sigma
  \quad \hbox{weakly star in $\H1\Vp\cap\L2V$}
  \non
  \\
  & \xieps \to \xi {\pier \ \hbox{ and } \Reps \to R}
  \quad \hbox{weakly in $\L2H$}
  \non
\end{align}
{\pier as $\eps \searrow 0,$} at least for a subsequence.
This implies, in particular, that the initial conditions 
for $(\mu,\phi,\sigma)$ are satisfied.
Moreover, the above convergence for $\phieps$
and the Aubin-Lions lemma (see, e.g., \cite[Thm.~5.1, p.~58]{Lions})
imply that {\pier
\Beq
\mueps \to \mu , \quad
\phieps \to \phi , \quad
\sigmaeps \to \sigma
  \quad \hbox{strongly in $\L2H$}.
\non 
\Eeq
}
Then, $\pi(\phieps)$ and $p(\phieps)$
converge to $\pi(\phi)$ and $p(\phi)$, respectively, strongly in $\L2H$.
Therefore, we can identify the limits of the nonlinear terms $\xieps$ and~$\Reps$.
For the former, we can apply, e.g., \cite[Cor.~2.4, p.~41]{Barbu}
{\pier and conclude that $\xi \in B(\varphi)$ a.e. in $Q.$}
For the latter we note that $\Reps$ converges to $p(\phi)(\sigma-\mu)$ 
{\pier strongly} in~$\LQ1${\pier ,
whence \eqref{defR} follows.}
At this point, we can write the integrated--in--time version of problem \accorpa{prima}{terza}
for the approximating solution
with time dependent test functions and take the limit as $\eps$ tends to zero.
We obtain the analogous systems for $(\mu,\phi,\sigma,\xi)$,
and this implies \accorpa{prima}{terza} for such a quadruplet.
This completes the proof of Theorem~\ref{Wellposedness}.
{\juerg \hfill{\qed}}


\section{Asymptotics}
\label{ASYMPTOTICS}
\setcounter{equation}{0}

This section is devoted to the proof of Theorems~\ref{Asymptotics} and~\ref{Error}.
In order to simplify the notation, 
we follow a general rule in performing our a~priori estimates.
The small-case italic $c$ without any subscript stands for different constants,
that may only depend on~$\Omega$, $T$, 
the shape of the nonlinearities 
and the norms of the initial data related to assumption~\eqref{hpdati}.
A~notation like~$c_\delta$ signals a constant that depends also on the parameter~$\delta$.
We point out that $c$ and $c_\delta$ do not depend on $\a$ and~$\b$
and that their meaning might change from line to line and even in the same chain of inequalities.
On the contrary, those constants we need to refer to are always denoted by different symbols,
e.g., by a capital letter.

\step
Proof of Theorem~\ref{Asymptotics}

We follow the argument done for \cite[{\pier Thm.}~2.6]{CGH} rather closely,
but we have to {\betti modify} the types of convergence
since our assumptions are different and more general.
We start {\gianni from \accorpa{genest}{genestbis}},
written for the solution $(\muab,\phiab,\sigmaab, {\betti \xiab})$, 
and improve the latter by estimating the norm of $\muab$ on its \rhs.
However, we omit the subscript{\betti s} $\a$ and $\b$ for a while. {\pier Thanks to 
\eqref{seconda} and \eqref{hp-pier}, we have that $|\xi |  \leq C \bigl( \Beta(\phi) + 1 \bigr)$ 
a.e. in $Q$. Then, by integrating over $\Omega$ we obtain}
{\gianni 
\Beq
  \iO |\xi(t)|
  \leq {}{\pier C}\! \iO \bigl( \Beta(\phi(t)) + 1 \bigr)
  \quad \aat.
  \label{stimaxi}
\Eeq
At this point, we can estimate the mean value of~$\mu$
{\pier on account of \eqref{hpPi} and \eqref{defWbeta}.  Indeed,
by just integrating \eqref{seconda} over~$\Omega$,
we deduce that}
\Bsist
  && \Bigl| \iO \mu(t) \Bigr|
  = \Bigl| \iO \bigl( \b\dt\phi +{\betti \xi+{\pier \pi}(\phi)} \bigr)(t) \Bigr|
  \non
  \\
  && \leq \b \norma{\dt\phi(t)}_1
  + c \bigl( \norma{\Beta(\phi(t))}_1 + \norma{\phi(t)}_1 + 1 \bigr) \vphantom\int
  \non
  \\
  && \leq {\pier c\,} \b^{1/2} \normaH{\dt\phi(t)}
  + c \bigl( \norma{F(\phi(t))}_1 + {\pier \normaH{\phi(t)}} + 1 \bigr) 
  \label{mediamu}
\Esist
{\pier \aat, beacause of the Lipschitz continuity of $\pi$ and the nonnegativity of $\Pi$.}}
Then, \eqref{genest} implies that the function
$t\mapsto\bigl|\textstyle\iO\mu(t)\,dt\bigr|$
is bounded in $L^2(0,T)$.
By combining this with~\eqref{genest} and the Poincar\'e inequality \eqref{poincare}, 
we derive that $\mu$ is bounded in~$\L2V$. {\betti
Hence, {\pier recalling estimates \eqref{genest}--\eqref{genestbis} it turns out that the convergences \accorpa{convmu}{pier1}
and a convergence for $\Ra $}} hold,
at least for a subsequence.
For the reader's convenience, we write {\juerg this conclusion explicitly},
as well as the {\betti consequences we are interested {\gianni in}{\pier. These are}
obtained by means of} strong compactness results
(see, e.g., \cite[Sect.~8, Cor.~4]{Simon}),
the Sobolev inequality \eqref{sobolev} and the \Lip\ continuity of $\pi$ and~$p$.
We~have {\pier
\begin{align}
  & \muab \to \mu
  \quad \hbox{weakly in $\L2V\cap\L2{{\pier \Lx6}}$}
  \non
  \\
  & \phiab \to \phi
  \quad \hbox{weakly star in $\L\infty V\cap\L2W$} \qquad
  \non
  \\
  & \sigmaab \to \sigma
  \quad  \hbox{weakly in $\H1\Vp\cap\L2V\cap\L2{{\pier \Lx6}}$}
  \non
  \\
  & \xiab \to \xi  \ \hbox{ and } \Reps \to R
  \quad  \hbox{weakly in $\L2H$}
  \non
  \\
  & \a\muab \to 0
  \quad \hbox{strongly in $\L\infty H {\pier {}\cap\L2V\cap\L2{{\pier \Lx6}}}$}
  \non
  \\
  & \b\dt\phiab \to 0
  \quad \hbox{strongly in $\L2H$}
  \non
  \\
  & \dt (\a\muab+\phiab) \to \dt\phi
  \quad \hbox{weakly in $\L2\Vp$}
  \non
   \\
  & \a\muab+\phiab \to \phi
  \quad \hbox{strongly in $\C0\Vp \cap \L2H$} 
  \non
  \\
  & \pi(\phiab) \to \pi(\phi) \ \hbox{ and } \ p(\phiab) \to p(\phi)
  \quad \hbox{strongly in $\L2H$}.
  \non
\end{align}
}
{\pier Hence, we} infer {\juerg that} $\phi$~and $\sigma$ satisfy the initial conditions~\eqref{cauchylim}.
Moreover, we deduce that {\betti $\xi\in B(\phi)$
(apply, e.g., \cite[Prop.~2.5, p.~27]{Brezis})}
and that $\Ra$ {\pier also converges to $p(\phi)(\sigma-\mu)$
weakly in~$\L1{\Lx p} $ for some $p\in (1,2)$: consequently, we have $R = p(\phi)(\sigma-\mu)$.}

Finally, we take the limit in the integrated--in--time version of problem \accorpa{prima}{terza}
for $(\muab,\phiab,\sigmaab,\xiab)$
with time-dependent test functions.
We obtain the analogue for the system \accorpa{primalim}{terzalim}.
Finally, as the solution of the limit problem is unique by Theorem~\ref{Uniqueness},
the convergence{\betti s} we have obtained for a subsequence
{\juerg hold} for the whole family.
This completes the proof of Theorem~\ref{Asymptotics}.\QED

\step
Proof of Theorem~\ref{Error}

As we use some ideas of~\cite{FGR},
it is convenient to rewrite the equations \eqref{prima} and \eqref{terza}
as abstract equations in the framework of the Hilbert triplet $(V,H,\Vp)$ 
related to an invertible operator.
To this end, we introduce the Riesz isomorphism $\A:V\to\Vp$
associated to the standard scalar product of~$V$, that~is
\Beq
  \< \A u, v > := (u,v)_V
  = \iO \bigl( \nabla u \cdot \nabla v + uv \bigr)
  \quad \hbox{for $u,v\in V$}.
  \label{defA}
\Eeq
We notice that $\A u=-\Delta u+u$ if $u\in W$
and that the restriction of $\A$ to $W$ is an isomorphism from $W$ onto~$H$.
We also remark that
\Bsist
  && \< \A u , \Am \vstar >
  = \< \vstar , u >
  \quad \hbox{for every $u\in V$ and $\vstar\in\Vp$}
  \label{propAuno}
  \\
  && \< \ustar , \Am \vstar >
  = (\ustar,\vstar)_*
  \quad \hbox{for every $\ustar,\vstar\in\Vp$,}
  \label{propAdue}
\Esist
where $(\cpto,\cpto)_*$ is the dual scalar product in $\Vp$
associated with the standard one in~$V$,
and recall that $\<\vstar,u>=\iO\vstar u$ if $\vstar\in H$.
As a consequence of~\eqref{propAdue}, we~have
\Beq
  \frac d{dt} \, {\betti \normaVp\vstar^2}
  = 2 \< \dt\vstar , \Am \vstar >
  \quad \hbox{for every $\vstar\in\H1\Vp$} .
  \label{propAtre}
\Eeq
{\pier In view of} the regularity conditions \Regsoluz\ and \Regsoluzlim, we rewrite 
\accorpa{prima}{terza} and \accorpa{primalim}{terzalim}
for the solution $(\muab,\phiab,\sigmaab)$ to \Pbl\ 
and the one of the limit problem, respectively.
If we term the latter $(\mul,\phil,\sigmal)$, we~have
\begin{align}
  & \a\dt\muab + \dt\phiab + \A\muab = \Ra + \muab
  \label{primaA}
  \\
  & \muab = \b\dt\phiab + \A\phiab + F'(\phiab) - \phiab
  \label{secondaA}
  \\
  & \dt\sigmaab + \A\sigmaab = -\Ra + \sigmaab
  \label{terzaA}
  \\
  & \dt\phil + \A\mul = \Rl + \mul
  \label{primalA}
  \\
  & \mul = \A\phil + F'(\phil) - \phil
  \label{secondalA}
  \\
  & \dt\sigmal + \A\sigmal = - \Rl + \sigmal ,
  \label{terzalA}
\end{align}
where $\Ra$ and $\Rl$ are defined by~\eqref{defR} according
to the equations we are considering.
All these equations are meant in $\Vp$ \aet.
However, \eqref{secondaA} and \eqref{secondalA}
also hold \aeQ.
Moreover, the solutions have to satisfy the initial conditions
\eqref{cauchy} and~\eqref{cauchylim}, respectively.
Now, we take the differences between 
\accorpa{primaA}{terzaA} and \accorpa{primalA}{terzalA}
at time $s\in(0,T)$
and test them~by
\Beq
  \Am (\a\muab + \phi)(s) , \quad
  -(\a\muab + \phi)(s),
  \aand
  \Am \sigma(s),
  \non
\Eeq
respectively, where we have set for convenience
\Beq
  \mu := \muab - \mul , \quad
  \phi := \phiab - \phil , \quad
  \sigma := \sigmaab - \sigmal,
  \aand
  R := \Ra - \Rl .
  \non
\Eeq
Next, we sum up and integrate over $(0,t)$ with respect to~$s$,
for an arbitrary $t\in(0,T)$.
We obtain
(by~omitting the evaluation at $s$ inside integrals, for brevity)
\Bsist
  && \iot \< \dt ( \a\muab + \phi ) , \Am (\a\muab + \phi) > \, ds
  + \iot \< \A\mu , \Am (\a\muab + \phi) > \, ds
  \non
  \\
  && \quad {}
  - \iot \< \mu , \a\muab + \phi > \, ds
  + \iot \< \b\dt\phiab , \a\muab + \phi > \, ds
  + \iot \< \A\phi , \a\muab + \phi > \, ds
  \non
  \\
  && \quad {}
  + \iot \< F'(\phiab) - F'(\phil) , \a\muab + \phi > \, ds
  \non
  \\
  && \quad {}
  + \iot \< \dt\sigma , \Am \sigma > \, ds
  + \iot \< \A\sigma , \Am \sigma > \, ds
  \non
  \\
  \separa
  && = \iot \< R , \Am (\a\muab + \phi) > \, ds
  + \iot \< \mu , \Am (\a\muab + \phi) > \, ds
  + \iot \< \phi , \a\muab + \phi > \, ds
  \non
  \\
  && \quad {}
  - \iot \< R , \Am \sigma > \, ds
  + \iot \< \sigma , \Am \sigma > \, ds \,.
  \non
\Esist
For the reader's convenience, we just rearrange and use the decomposition $F'=B+\pi$.
We have
\Bsist
  && \iot \< \dt \bigl( \a\muab + \phi \bigr) , \Am (\a\muab + \phi) > \, ds
  \non
  \\
  && \quad {}
  + \iot \< \A\mu , \Am (\a\muab + \phi) > \, ds
  - \iot \< \mu , \a\muab + \phi > \, ds
  \non
  \\
  && \quad {}
  + \iot \< \A\phi , \phi > \, ds
  + \iot \< B(\phiab) - B(\phil) , \phi > \, ds
  \non
  \\
  && \quad {}
  + \iot \< \dt\sigma , \Am \sigma > \, ds
  + \iot \< \A\sigma , \Am \sigma > \, ds
  \non
  \\
  \separa
  && =
  \iot \< R , \Am (\a\muab + \phi - \sigma) > \, ds
  - \iot \< \b\dt\phiab , \a\muab + \phi > \, ds
  - \iot \< \A\phi , \a\muab > \, ds
  \non
  \\
  && \quad {}
  + \iot \< \mu , \Am (\a\muab + \phi) > \, ds
  - \iot \< \pi(\phiab) - \pi(\phil) , \phi > \, ds
  \non
  \\
  && \quad {}
  - \iot \< F'(\phiab) - F'(\phil) , \a\muab > \, ds
  + \iot \< \phi , \a\muab + \phi > \, ds
  + \iot \< \sigma , \Am\sigma > \, ds.
  \non
\Esist
At this point, we account for \accorpa{defA}{propAtre}
and observe that the second and third terms on the \lhs\ cancel out.
Finally, {\betti owing {\juerg to} the initial conditions
$(\a\muab+\phi)(0)=\a\muz$ and $\sigma(0)=0$, we deduce}
\Bsist
  && \frac 12 \, \normaVp{(\a\muab + \phi)(t)}^2
  + \iot \normaV\phi^2 \, ds
  + \intQt \bigl( B(\phiab) - B(\phil) \bigr) \phi 
  \quad
  \non
  \\
  && \quad {}
  + \frac 12 \, \normaVp{\sigma(t)}^2 
  + \intQt |\sigma|^2
  \non
  \\ 
  && = \frac 12 \, \normaVp{\a\muz}^2
  + \iot \bigl( R , \a\muab + \phi - \sigma \bigr)_* \, ds  
  - \iot \< \b\dt\phiab , \a\muab + \phi > \, ds             
  \non
  \\
  && \quad {}
  - \iot \bigl( \phi , \a\muab \bigr)_V \, ds               
  + \iot \bigl( \mu , \a\muab + \phi \bigr)_* \, ds         
  - \intQt \bigl( \pi(\phiab) - \pi(\phil) \bigr) \, \phi   
  \non
  \\
  && \quad {}
  - \intQt \bigl( F'(\phiab) - F'(\phil) \bigr) \, \a\muab   
  + \intQt \phi (\a\muab + \phi)                            
  + \iot \normaVp\sigma^2 \, ds \,.                        
  \label{pererrore}
\Esist
All {\juerg of} the terms on the \lhs\ are nonnegative,
the third one by monotonicity.
Now, we treat each integral on the \rhs, separately.
In the sequel, $\delta$~is a positive parameter 
whose value will be chosen at the end of the procedure.
We first observe that \eqref{genest} holds for the solution $(\muab,\phiab,\sigmaab)$
and that Theorem~\ref{Asymptotics} improves~\eqref{genestbis} for such a solution.
Indeed, the restricted setting of regular potentials
satisfying \eqref{classico} led to~\eqref{mediamu}.
So, as we have seen in the previous proof, {\pier \eqref{genest} and \eqref{genestbis}~imply}
\Beq
  \norma\muab_{\L2V}
  + \norma\phiab_{\L2W}
  \leq c.
  \label{genestter}
\Eeq
Now, we prepare estimates for $\normaVp\mu$ and $\normaVp R$ \aet.
Again for simplicity, in performing them,
we omit writing the evaluation point.
From the mean value theorem and the third assumption {\pier in}~\eqref{classico}
we easily derive that
\Beq
  |F'(\phiab)-F'(\phil)|
  \leq c|\phi|(|\phiab|^4+|\phil|^4+1)
  \quad \aeQ.
  \non
\Eeq
Therefore, by the \holder\ and Sobolev inequalities, we infer~that
\begin{align}
  & \norma{F'(\phiab)-F'(\phil)}_{6/5}
  \leq c \, \norma\phi_6 \bigl( \norma{\phiab^4}_{3/2} + \norma{\phil\,^4}_{3/2} + 1 \bigr)
  \non
  \\
  & = c \, \norma\phi_6 \bigl( \norma\phiab_6^4 + \norma\phil_6^4 + 1 \bigr)
  \leq c \normaV\phi \bigl( \normaV\phiab^4 + \normaV\phil^4 + 1 \bigr)
  \leq c \normaV\phi \, ,
  \label{stimadiffF}
\end{align}
the last inequality {\betti following from} estimate \eqref{genest} for $\phiab$ and the regularity \eqref{regphilim} of~$\phil$.
{\betti Taking the} difference between \eqref{secondaA} and \eqref{secondalA}
and {\betti using} the dual Sobolev inequality~\eqref{dualsobolev},
we deduce~that
\Bsist
  && \normaVp\mu
  = \normaVp{\b\dt\phiab + \A\phi + F'(\phiab) - F'(\phil) - \phi}
  \non
  \\
  && \leq \b \normaVp{\dt\phiab} + \normaV\phi + c \norma{F'(\phiab) - F'(\phil)}_{6/5} + \normaVp\phi 
  \non
  \\
  && \leq c \bigl(
    \b \normaVp{\dt\phiab} + \normaV\phi
  \bigr) 
  \leq c \, \b^{1/2} + c \, \normaV\phi \, , 
  \label{stimamu}
\Esist
{\juerg where the last inequality follows from} \eqref{genest}.
In order to estimate~$\normaVp R$, we first observe that
the boundedness and the \Lip\ continuity of~$p$
and the Sobolev inequality (applied to $\nabla\phil$ and~{\betti the test function $v\in V$})
imply {\juerg that, for every $v\in V$,}
\Bsist
  && \normaV{p(\phil)v}
  \leq \normaH{p(\phil)v}
  + \normaH{\nabla p(\phil) \, v}
  + \normaH{p(\phil)\nabla v}
  \non
  \\
  && \leq c \normaH v + c \norma{\nabla\phil}_4 \, \norma v_4 + c \normaH{\nabla v}
  \leq c (\normaW\phil \, + 1) \normaV v \,.
  \non
\Esist
Hence, we have for every $v\in V$ {\juerg the estimate} 
\Bsist
  && \iO R v
  = \iO \bigl( p(\phiab)(\sigmaab-\muab) - p(\phil)(\sigmal-\mul) \bigr) v
  \non
  \\
  && \leq \iO |p(\phiab) - p(\phil)| \, |\sigmaab-\muab| \, |v|
  + \Bigl| \iO p(\phil) (\sigma-\mu) \, v \Bigr|
  \non
  \\
  && \leq c \norma\phi_3 \, \norma{\sigmaab-\muab}_3 \norma v_3
  + \normaVp{\sigma-\mu} \, \normaV{p(\phil)\,v}
  \non
  \\
  && \leq c \normaV\phi \, \normaV{\sigmaab-\muab} {\betti  \normaV v}
    + c \normaVp{\sigma-\mu} \, (\normaW\phil + 1) \normaV v \,.
  \non
\Esist
Therefore, we can estimate $\normaVp R$ \aet {\betti , also} owing to \eqref{stimamu}, in  this~way:
\Bsist
  && \normaVp R
  \leq c \normaV\phi \, \normaV{\sigmaab-\muab} 
  + c \normaVp{\sigma-\mu} \, (\normaW\phil + 1)
  \non
  \\
  && \leq c \normaV\phi \, \normaV{\sigmaab-\muab} 
  + c \bigl( \normaVp\sigma + \b^{1/2} + \normaV\phi \bigr) (\normaW\phil + 1)
  \non
  \\
  && \leq c \, \psia (\normaV\phi + \normaVp\sigma)
  + c \, \b^{1/2} \, \psil \,,
  \non
\Esist
where $\psia,\psil:(0,T)\to\erre$ are defined~by
\Beq
  \psia := \normaV{\sigmaab-\muab} + \normaW\phil + 1
  \aand
  \psil := \normaW\phil + 1,
  \quad \aet \,.
  \non
\Eeq
Coming back to the \rhs\ of~\eqref{pererrore},
we can treat the first term as follows:
\Bsist
  && \iot \bigl( R , \a\muab + \phi - \sigma \bigr)_* \, ds
  \leq \iot \normaVp R \, \normaVp{\a\muab + \phi - \sigma} \, ds
  \non
  \\
  && \leq \iot \bigl(
    c \psia (\normaV\phi + \normaVp\sigma)
    + c \b^{1/2} \, \psil
  \bigr)
  \bigl(
    \normaVp{\a\muab+\phi}
    + \normaVp\sigma
  \bigr) \, ds
  \non
  \\
  && \leq \delta \iot \normaV\phi^2 \, ds
  + \b \ioT |\psil|^2 \, ds
  + c_\delta \iot \psia^2 \bigl( \normaVp{\a\muab+\phi}^2 + \normaVp\sigma^2 \bigr) \, ds \,.
  \qquad
  \label{stimaf}
\Esist
We observe at once that the regularity {\pier \eqref{regphilim} for $\phil$}
and estimates \eqref{genest} for $\sigmaab$ and \eqref{genestter} for $\muab$ imply that 
$\psil\in L^2(0,T)$ and that
$\psia$ is bounded in $L^2(0,T)$, {\juerg so that} $\psia^2$ is bounded in $L^1(0,T)$.
This will allow us to apply the Gronwall lemma.
Now, we estimate the next term on the \rhs\ of \eqref{pererrore}.
{\pier Using \eqref{genest}, we see that}
\Bsist
  && - \iot \< \b\dt\phiab , \a\muab + \phi > \, ds
  \non
  \\
  && \leq \a^2 \, \norma\muab_{\L2H}^2 
  + \delta \iot \normaH\phi^2 \, ds
  + c_\delta \, \b^2 \, \norma{\dt\phiab}_{\L2H}^2 
  \non
  \\
  && \leq c \a
  + \delta \iot \normaV\phi^2 \, ds
  + c_\delta \, \b \,.
  \label{stimag}
\Esist
Next, we have
\Beq
  - \iot \bigl( \phi , \a\muab \bigr)_V \, ds
  \leq \delta \iot \normaV\phi^2 \, ds
  + c_\delta \, \a^{{\betti 2}}
  \label{stimad}
\Eeq
thanks to \eqref{genestter} for $\muab$,
as well as, by \eqref{stimamu},
\Bsist
  && \iot \bigl( \mu , \a\muab + \phi \bigr)_* \, ds
  \leq c \iot \bigl( \b^{1/2} + \normaV\phi \bigr) \normaVp{\a\muab + \phi} \, ds
  \non
  \\
  && \leq \delta \iot \normaV\phi^2
  + \b
  + c_\delta \iot \normaVp{\a\muab + \phi}^2 \, ds \,.
  \label{stimae}
\Esist
Moreover, by using the \Lip\ continuity of~$\pi$,
the interpolation inequality~\eqref{interpol}
and \eqref{genestter} for~$\muab$ once more,
we can write
\Bsist
  && - \intQt \bigl( \pi(\phiab) - \pi(\phil) \bigr) \, \phi
  \leq c \intQt |\phi|^2
  \leq c \iot \normaV\phi \, \normaVp\phi \, ds
  \non
  \\
  && \leq \delta \iot \normaV\phi^2 \, ds
  + c_\delta \iot \normaVp\phi^2 \, ds
  \non
  \\
  && \leq \delta \iot \normaV\phi^2 \, ds
  + {\betti c_\delta \iot \normaVp{\a\muab+\phi}^2 \, ds} {\pier{} + c_\delta \a^2} \, .
  \label{stimac}
\Esist
The next term to {\juerg deal} with is the one involving~$F'$.
We use \eqref{stimadiffF}, the \holder, Sobolev and Young inequalities,
and {\pier the} estimate \eqref{genestter} for~$\muab${\pier . Thus, we} have
\Bsist
  && - \intQt \bigl( F'(\phiab) - F'(\phil) \bigr) \, \a\muab 
  \leq \iot \norma{F'(\phiab)-F'(\phil)}_{6/5} \, \norma{\a\muab}_6 \, ds
  \non
  \\
  && \leq c \iot \normaV\phi \, \a \normaV\muab \, ds
  \leq \delta \iot \normaV\phi^2 \, ds
  + c_\delta \, \a^{{\betti 2}} \,.
  \label{stimam}
\Esist
Finally, the last integral on the \rhs\ of \eqref{pererrore} does not need any treatment 
and the {\juerg preceding} term can be estimated in this way:
\Beq
  \iot \< \phi , \a\muab + \phi > \, ds
  \leq \delta \iot \normaV\phi^2 \, ds
  + c_\delta \iot \normaVp{\a\muab + \phi}^2 \, ds .
  \label{stimab}
\Eeq
At this point, we combine \eqref{pererrore}
and the list \accorpa{stimaf}{stimab} of estimates we have obtained.
Then, we choose $\delta$ small enough, recall that
$\psil\in L^2(0,T)$ and that $\psia^2$ is bounded in~$L^1(0,T)$,
and apply the Gronwall lemma in the form \cite[Lemma~A.4, p.~156]{Brezis}.
We obtain
\Beq
  \frac 12 \, \normaVp{(\a\muab + \phi)(t)}^2
  + \iot \normaV\phi^2 \, ds
  + \frac 12 \, \normaVp{\sigma(t)}^2 
  + \intQt |\sigma|^2
  \leq c (\a+\b) 
  \non
\Eeq
for every $t\in[0,T]$.
As $\normaVp{\a\muab(t)}^2\leq c\a$ for every $t\in[0,T]$ by \eqref{genest},
the above inequality implies 
\begin{align}
  \norma{\phi}_{\L\infty\Vp\cap\L2V}
  + \norma{\sigma}_{\L\infty\Vp\cap\L2H}
  \leq C \, \bigl( \a^{1/2} + \b^{1/2} \bigr). 
  \label{pier2}
\end{align}
{\pier Now, we take the differences of equations \eqref{secondaA} and \eqref{secondalA}
and estimate the $\L2\Vp $ norm of it. With the help of \eqref{genest} and \eqref{stimadiffF}
it is straightforward to infer that 
\begin{align}
  &\norma{\mu}_{\L2\Vp} \leq c\, \b \, \norma{\dt\phiab}_{\L2H} + 
   c \,  \norma{\phi}_{\L2V} 
   \non
   \\
 &\quad  {}+  c \, \norma{F'(\phiab)-F'(\phil)}_{\L2{L^{6/5} (\Omega)}} + \norma{\phi}_{\L2\Vp} 
  \non
   \\
 &{}\leq   c \, \b^{1/2}  +    c \,  \norma{\phi}_{\L2V}  .
  \label{pier3} 
\end{align} 
Hence, in view of \eqref{pier2} and \eqref{pier3} we finally obtain the estimate
\eqref{error}, where one has to {\juerg read} $\phil, \, \mul$ and $\sigmal$
in place of $\phi, \, \mu $} and~$\sigma$, respectively,
due to the change of notations within this proof.\QED

\Brem
\label{Uniqueness}
By going through the above proof, one sees that uniqueness for the limit problem
\accorpa{primalim}{cauchylim} has been never used,
that is, the following formulation of Theorem~\ref{Error} has been proved:
the error estimate \eqref{error} holds for every solution $(\mu,\phi,\sigma)$
of the limit problem satisfying the regularity requirements \Regsoluzlim.
This implies the uniqueness for such a solution.
Indeed, if $(\mu_i,\phi_i,\sigma_i)$, $i=1,2$, are two solutions {\betti of the limit problem},
by writing \eqref{error} for both of them 
and using uniqueness for the solution $(\phiab,\muab,\sigmaab)$, 
one immediately derives
\begin{align}
  &\norma{\phi_1-\phi_2}_{\L\infty\Vp\cap\L2V}
 {\pier{} +  \norma{\mu_1-\mu_2}_{\L2\Vp}}
  \non
  \\
   &\quad + \norma{\sigma_1 - \sigma_2}_{\L\infty\Vp\cap\L2H}
  \leq C \, \bigl( \a^{1/2} + \b^{1/2} \bigr) 
  \non 
\end{align}
for every $\a,\b\in(0,1)$,
whence $\phi_1=\phi_2$,  {\pier $\mu_1=\mu_2$} and $\sigma_1=\sigma_2$.
Then, by comparison in \eqref{secondalim}, it follows that {\pier $\xi_1=\xi_2$}, as well.
\Erem



\vspace{3truemm}

\Begin{thebibliography}{10}

\bibitem{Barbu}
{\sc V. Barbu},
``Nonlinear Differential Equations of Monotone Types in Banach spaces'',
Springer Monographs in Mathematics, 2010.

\bibitem{Brezis}
{\sc H. Brezis,}
``Op\'erateurs maximaux monotones et semi-groupes de contractions
dans les espaces de Hilbert'',
North-Holland Math. Stud.~{\bf 5},
North-Holland, Amsterdam, 1973.

\bibitem{CH} 
{\sc J.W.~Cahn and J.E.~Hilliard}, 
{\em Free energy of a nonuniform system I. Interfacial free energy}, 
J. Chem. Phys., {\bf 2} (1958), pp.~258--267.

{\betti 
\bibitem{CBCB} 
{\sc C. Chatelain, T. Balois, P. Ciarletta, and M. Ben Amar}, 
{\em Emergence of microstructural patterns in skin cancer: a phase separation
analysis in a binary mixture}, New J. Phys., \textbf{13} (2011), 115013 (21 pp.).
}

\bibitem{CGH}
{\sc P. Colli, G. Gilardi, {\betti and } D. Hilhorst},
{\em On a Cahn--Hilliard type phase field system related to tumor growth},
{\gianni{\pier Discrete Contin. Dyn. Syst.,} \textbf{35} (2015), pp.~2423--2442.}

{\betti 
\bibitem{CLLW} 
{\sc V. Cristini, X. Li, J.S. Lowengrub, and S.M. Wise}, 
{\em Nonlinear simulations of solid tumor growth using a mixture model: invasion and branching},
 J. Math. Biol., \textbf{58}  (2009), pp.~723--763.

\bibitem{CL} 
{\sc V. Cristini and J. Lowengrub}, ``Multiscale modeling of cancer. An Integrated Experimental and Mathematical Modeling Approach'',
Cambridge Univ. Press, Cambridge, 2010.
}

{\pier
\bibitem{DHK}
{\sc C. Dupaix, D. Hilhorst  and I.N. Kostin,} 
{\em The viscous Cahn-Hilliard equation as a limit
of the phase field model: lower semicontinuity of the attractor}, 
J. Dynam. Differential Equations, \textbf{11}  (1999), pp.~333-353.

\bibitem{DHL}
{\sc C. Dupaix, D. Hilhorst  and P.  Lauren\c cot,} 
{\em Upper-semicontinuity of the attractor for a
singularly perturbed phase field model}, 
Adv. Math. Sci. Appl., \textbf{8} (1998), pp.~115-143.
}

\bibitem{EllSt} 
{\sc C.M. Elliott and A.M. Stuart}, 
{\em Viscous Cahn--Hilliard equation. II. Analysis}, 
J. Differential Equations,
{\bf 128} (1996),  pp.~387--414.

\bibitem{EllZh} 
{\sc C.M. Elliott and S. Zheng}, 
{\em On the Cahn--Hilliard equation}, 
Arch. Rational Mech. Anal.,
{\bf 96} (1986), pp.~339--357.

{\betti 
\bibitem{Fetal} 
{\sc H.B. Frieboes, F. Jin, Y.-L. Chuang, S.M. Wise, J.S. Lowengrub, and V. Cristini}, 
{\em Three-dimensional multispecies nonlinear tumor growth-II: Tumor invasion and angiogenesis}, 
J. Theoret. Biol., \textbf{264}  (2010), pp.~1254--1278.
}

\bibitem{FGR}
{\sc S. Frigeri, M. Grasselli, and E. Rocca},
{\em On a diffuse interface model of tumor growth},
{\bet European J. Appl. Math., DOI: 10.1017/S0956792514000436}

{\betti 
\bibitem{HZO} 
{\sc A. Hawkins-Daarud, K.G. van der Zee, and J.T. Oden}, 
{\em Numerical simulation
of a thermodynamically consistent four-species tumor growth model}, 
Int. J. Numer. Meth. Biomed. Engng., \textbf{28} (2011), pp.~3--24.

\bibitem{HPZO} 
{\sc A. Hawkins-Daarud, S. Prudhomme, K.G. van der Zee, and J.T. Oden}, 
{\em Bayesian calibration, validation, and uncertainty quantification of diffuse interface models of tumor growth}, 
J. Math. Biol., \textbf{67}  (2013), pp.~1457--1485.

\bibitem{HKNV} 
{\sc D. Hilhorst, J. Kampmann, T.N. Nguyen, and K.G. Van der Zee}, {\bet {\em Formal
asymptotic limit of a diffuse-interface tumor-growth model}},
Math. Models Methods Appl. Sci., {\bet DOI: 10.1142/S0218202515500268}
}

\bibitem{Lions}
{\sc J.-L.~Lions},
``Quelques m\'ethodes de r\'esolution des probl\`emes
aux limites non lin\'eaires'',
Dunod; Gauthier-Villars, Paris, 1969.

{\betti 
\bibitem{Letal} 
{\sc J.S. Lowengrub, H.B. Frieboes, F. Jin, Y.-L. Chuang, X. Li, P. Macklin, S.M. Wise, and V. Cristini}, 
{\em Nonlinear modelling of cancer: bridging the gap between cells and tumours}, Nonlinearity, \textbf{23}  (2010), pp.~R1--R91.

\bibitem{LTZ} 
{\sc J. Lowengrub, E. Titi, and K. Zhao}, 
{\em Analysis of a mixture model of tumor growth}, 
European J. Appl. Math., \textbf{24} (2013), pp.~1--44.

\bibitem{OHP} 
{\sc J.T. Oden, A. Hawkins, and S. Prudhomme}, 
{\em General diffuse-interface theories and an approach to predictive tumor growth modeling}, 
Math. Models Methods Appl. Sci., 
\textbf{20} (2010), pp.~477--517.

\bibitem{OPH} 
{\sc J.T. Oden, E.E. Prudencio, and A. Hawkins-Daarud}, 
{\em Selection and assessment of phenomenological models of tumor growth},
Math. Models Methods Appl. Sci., \textbf{23}  (2013), pp.~1309--1338.
}

{\pier
\bibitem{Ros1}
{\sc R. Rossi},
{\em Asymptotic analysis of the Caginalp phase-field model 
for two vanishing time relaxation parameters}, 
Adv. Math. Sci. Appl., \textbf{13} (2003), pp.~249--271. 

\bibitem{Ros2}
{\sc R. Rossi},
{\em On two classes of generalized viscous Cahn-Hilliard equations},  
Commun. Pure Appl. Anal., \textbf{4} (2005), pp.~405--430. 
}

\bibitem{Simon}
{\sc J. Simon},
{\em Compact sets in the space $L^p(0,T; B)$},
Ann. Mat. Pura Appl.~(4), {\bf 146} (1987), pp.~65-96.

{\betti
\bibitem{WZ} 
{\sc X. Wang and Z. Zhang}, 
{\em  Well-posedness of the Hele--Shaw--Cahn--Hilliard system}, Ann. Inst. H. Poincar\'{e} Anal. Non Lin\'{e}aire,
\textbf{30} (2013), pp.~367--384.

\bibitem{WLFC} 
{\sc S.M. Wise, J.S. Lowengrub, H.B. Frieboes, and V. Cristini}, 
{\em Three-dimensional multispecies nonlinear tumor growth-I: Model and numerical method}, 
J. Theoret. Biol., \textbf{253}  (2008),  pp.~524--543.
 
\bibitem{WZZ} 
{\sc X. Wu, G.J. van Zwieten, and K.G. van der Zee}, 
{\em Stabilized second-order convex splitting schemes for Cahn--Hilliard
models with applications to diffuse-interface tumor-growth models}, 
Int. J. Numer. Meth. Biomed. Engng., \textbf{30} (2014), pp.~180--203.
}

\End{thebibliography}

\End{document}

\bye